\documentclass[a4paper]{amsart}
\usepackage[english]{babel}

\usepackage{amsmath,amssymb,amsthm}
\theoremstyle{plain}
\newtheorem{theorem}{Theorem}[section]
\newtheorem{corollary}[theorem]{Corollary}
\newtheorem{prop}[theorem]{Proposition}
\newtheorem{lemma}[theorem]{Lemma}

\theoremstyle{definition}
\newtheorem{remark}[theorem]{Remark}
\newtheorem{example}[theorem]{Example}

\newcommand{\ecc}{\ensuremath{\overline{\mathrm{co}}\,}}
\newcommand{\ext}[1][X^*]{\ensuremath{\mathrm{ex}\left(B_{#1}\right)}}
\newcommand{\extr}{\ensuremath{\mathrm{ex}}}

\DeclareMathOperator{\re}{Re\,}

\newcounter{equi1}
\newenvironment{equi}
{\begin{list}
        {(\roman{equi1})}
        {\setlength{\itemsep}{0ex plus 0.2ex minus 0ex}
         \setlength{\topsep}{0ex}
         \setlength{\parsep}{0ex}
         \setlength{\labelwidth}{7ex}
         \usecounter{equi1}}}
{\end{list}}

\begin{document}

\thispagestyle{empty}

\title[Daugavet property of $C^*$-algebras and $JB^*$-triples]{The
Daugavet Property of $\mathbf{C^*}$-algebras,
$\mathbf{JB^*}$-triples, and of their isometric preduals}
\date{November 18th, 2004}
\maketitle

\thispagestyle{empty}

\bigskip

\bigskip

 \centerline{\textsc{\large Julio Becerra Guerrero}
 \footnote{Partially supported by Junta de Andaluc\'{\i}a grant
 FQM-0199}}

\begin{center} Departamento de Matem\'{a}tica Aplicada \\ Facultad de
Ciencias \\ Universidad de Granada \\ 18071 Granada, SPAIN\\
\emph{Email address:} \texttt{juliobg@ugr.es} \end{center}

\bigskip

\bigskip

\bigskip

 \centerline{\textsc{\large Miguel Mart\'{\i}n}
\footnote{Partially supported by Spanish MCYT project no.\
BFM2003-01681} \footnote{Corresponding author. \emph{Email:}
\texttt{mmartins@ugr.es}, \emph{Fax:} \texttt{+34958243272}}}

\begin{center} Departamento de An\'{a}lisis Matem\'{a}tico \\ Facultad de
Ciencias \\ Universidad de Granada \\ 18071 Granada, SPAIN \\
\emph{Email address:} \texttt{mmartins@ugr.es} \end{center}

\vspace{2.5cm}

{\Large \textbf{Abstract:}}

\noindent A Banach space $X$ is said to have the Daugavet property
if every rank-one operator $T:X\longrightarrow X$ satisfies $\|Id
+ T\| = 1 + \|T\|$. We give geometric characterizations of this
property in the settings of \mbox{$C^*$-algebras}, $JB^*$-triples
and their isometric preduals. We also show that, in these
settings, the Daugavet property passes to ultrapowers, and thus,
it is equivalent to an stronger property called the uniform
Daugavet property.

\bigskip

\bigskip

{\large \textbf{Keywords:} \footnotesize

\noindent $C^*$-algebra; von Neumann predual; $JB^*$-triple;
predual of a $JBW^*$-triple; Daugavet equation; Daugavet property;
rough norm; Fr\'{e}chet-differentiability.


\markboth{The Daugavet property of $C^*$-algebras and
$JB^*$-triples}{J.~Becerra and M.~Mart\'{\i}n}

\newpage
\normalsize

\section{Introduction}
A Banach space $X$ is said to have the \emph{Daugavet property}
\cite{KSSW} if every rank-one operator $T:X\longrightarrow X$
satisfies the norm identity
\begin{equation*}\label{DE}\tag{\textrm{DE}}
\|Id + T\| = 1 + \|T\|,
\end{equation*}
known as \emph{Daugavet equation}. In such a case, all weakly
compact operators on $X$ also satisfy (\ref{DE}) (see
\cite[Theorem~2.3]{KSSW}). Therefore, this definition of Daugavet
property coincides with those that appeared in \cite{Chau} and
\cite{AAB}.

The study of the Daugavet equation was inaugurated by I.~Daugavet
\cite{Dau} in 1961 by proving that every compact operator on
$C[0,1]$ satisfies (\ref{DE}). Over the years, the validity of the
Daugavet equation was proved for compact operators on various
spaces, including $C(K)$ and $L_1(\mu)$ provided that $K$ is
perfect and $\mu$ does not have any atoms (see \cite{Wer0} for an
elementary approach), and certain function algebras such as the
disk algebra $A(\mathbb{D})$ or the algebra of bounded analytic
functions $H^\infty$ \cite{WerJFA,Woj}. In the nineties, new ideas
were infused into the field and the geometry of Banach spaces
having Daugavet property was studied. The state-of-the-art on the
subject can be found in \cite{KSSW,WerSur}. For very recent
results we refer the reader to \cite{BKSW,KaKaWe,KaWe} and
references therein.

Let us mention here several facts concerning the Daugavet property
which are relevant to our discussion. It is clear that $X$ has the
Daugavet property whenever its topological dual $X^*$ does, but
the converse result is false ($X=C[0,1]$, for instance). It is
known that a space with the Daugavet property cannot have the
Radon-Nikod\'{y}m property (RNP in short) \cite{Woj}; even more, every
weakly open subset of its unit ball has diameter 2 \cite{Shv}. A
space with the Daugavet property contains a copy of $\ell_1$
\cite{KSSW}, it does not have an unconditional basis \cite{Kad}
and it does not even embed into a space with an unconditional
basis \cite{KSSW}.

In 2002, T.~Oikhberg \cite{Oik} carried the classical results on
the Daugavet property for $C(K)$ and $L_1(\mu)$ to the
non-commutative case, characterizing when (complex) $C^*$-algebras
and preduals of von Neumann algebras have the Daugavet property. A
$C^*$-algebra has the Daugavet property if and only if it does not
have atomic projections; if the algebra is a von Neumann algebra
(i.e., it is a dual space), its (unique) isometric predual has the
Daugavet property if and only if the algebra does. In 2004,
T.~Oikhberg and the second named author \cite{MaOi}, translated
these results to the non-associative case, characterizing
(complex) $JB^*$-triples and predual of (complex) $JBW^*$-triples
having the Daugavet property in an analogous way, replacing atomic
projections by minimal tripotents. The necessary definitions and
basic results on $JB^*$-triples are presented in
section~\ref{sec-JB}.

In the present paper we give geometric characterizations of the
Daugavet property in the setting of real and complex
$JB^*$-triples and their isometric preduals. In particular, our
results contain the already mentioned ones of \cite{MaOi,Oik} for
complex $C^*$-algebras and complex $JB^*$-triples, but our proofs
are independent.

To state the main results of the paper we need to fix notation and
recall some definitions.

Let $X$ be a Banach space. The symbols $B_X$ and $S_X$ denote,
respectively, the closed unit ball and the unit sphere of $X$. Let
us fix $u$ in $S_X$. We define the set $D(X,u)$ of all
\emph{states} of $X$ relative to $u$ by $$D(X,u):= \{f\in B_{X^*}\
: \ f(u)=1\},$$ which is a non-empty $w^*$-closed face of
$B_{X^*}$. The norm of $X$ is said to be \emph{smooth} at $u$ if
$D(X,u)$ reduces to a singleton, and it is said to be
\emph{Fr\'{e}chet-smooth} or \emph{Fr\'{e}chet differentiable} at $u\in
S_X$ whenever there exists $\displaystyle \lim _{\alpha
\rightarrow 0} \frac{\Vert u+\alpha x\Vert -1}{\alpha}$ uniformly
for $x\in B_X$. We define the \emph{roughness of $X$ at $u$} by
the equality
$$ \eta (X,u)~:=~\limsup_{\Vert h\Vert \rightarrow 0}
\dfrac{\Vert u+h\Vert +\Vert u-h\Vert -2}{\Vert h\Vert }.$$ We
remark that the absence of roughness of $X$ at $u$ (i.e., $\eta
(X,u)=0$) is nothing other than the Fr\'{e}chet-smoothness of the norm
of $X$ at $u$ \cite[Lemma~I.1.3]{DGZ}. Given $\delta >0$, the
Banach space $X$ is said to be \emph{$\delta$-rough} if, for every
$u$ in $S_X$, we have $\eta (X,u)\geqslant \delta $. We say that
$X$ is \emph{rough} whenever it is $\delta$-rough for some $\delta
>0$, and \emph{extremely rough} whenever it is $2$-rough. Roughly
speaking, the space $X$ is rough if its norm is ``uniformly''
non-differentiable at any point. A \emph{slice} of $B_X$ is a
subset of the form $$S(B_X,f,\alpha)=\bigl\{x\in B_X\ : \ \re f(x)
> 1 - \alpha \bigr\},$$ where $f\in S_{X^*}$ and $0<\alpha<1$. If
$X$ is a dual space and $f$ is actually taken from the predual, we
say that $S(B_X,f,\alpha)$ is a $w^*$-slice. By
\cite[Proposition~I.1.11]{DGZ}, the norm of $X$ is $\delta$-rough
if and only if every nonempty $w^*$-slice of $B_{X^*}$ has
diameter greater or equal than $\delta$.

Finally, a point $x\in S_X$ is said to be an \emph{strongly
exposed point} if there exists $f\in D(X,x)$ such that $\lim \|x_n
- x\|=0$ for every sequence $(x_n)$ of elements of $B_{X}$ such
that $\lim\,\re f(x_n)=1$ (equivalently, there are slices defined
by $f$ with arbitrary small diameter). It is known that $x$ is
strongly exposed if and only if there is a point of
Fr\'{e}chet-smoothness in $D(X,x)$ (see \cite[Corollary~I.1.5]{DGZ}).

The main results of the paper are the characterizations of the
Daugavet property for $JB^*$-triples and preduals of
$JBW^*$-triples given in Theorems \ref{th-JB} and \ref{th-jbw}
respectively. For a real or complex $JB^*$-triple $X$, the
following are equivalent:
\begin{enumerate}
\item[$(i)$] $X$ has the Daugavet property,
\item[$(ii)$] the norm of $X$ is extremely rough,
\item[$(iii)$] the norm of $X$ is not Fr\'{e}chet-smooth at any point.
\end{enumerate}
For the predual $X_*$ of a real or complex $JBW^*$-triple $X$, the
following are equivalent:
\begin{enumerate}
\item[$(i)$] $X$ has the Daugavet property,
\item[$(ii)$] $X_*$ has the Daugavet property,
\item[$(iii)$] every relative weak-open subset of $B_{X_*}$ has diameter $2$,
\item[$(iv)$] $B_{X_*}$ has no strongly exposed points,
\item[$(v)$] $B_{X_*}$ has no extreme points.
\end{enumerate}
This characterizations allow us to prove that, for $JB^*$-triples
and for preduals of $JBW^*$-triples, the Daugavet property passes
to ultrapowers. As a consequence, a stronger version of the
Daugavet property introduced in \cite{BKSW}, called the uniform
Daugavet property, is equivalent to the usual Daugavet property in
the setting of $JB^*$-triples and their isometric preduals.

The outline of the paper is as follows. In section~2 we give
sufficient conditions for a Banach space to have the Daugavet
property, which will be the keys to state the rest of the paper.

Section~3 is devoted to the above cited characterizations of the
Daugavet property for real or complex $JB^*$-triples and their
isometric preduals, and we dedicate section~4 to particularize
these result to the setting of real or complex $C^*$-algebras and
von Neumann preduals.

Finally, in section~5 we study the behaviour of the Daugavet
property for ultraproducts of $JB^*$-triples and of preduals of
$JBW^*$-triples. As a consequence, we show that the already
mentioned uniform Daugavet property and the Daugavet property
coincide in real or complex $JB^*$-triples and their isometric
preduals.

Throughout the paper, for a subset $A$ of a Banach space, we write
$\ecc(A)$ for the closed convex hull of $A$, we use $\extr(B)$ to
denote the set of extreme points of the convex set $B$ and,
finally, if $X$ and $Y$ are Banach spaces, we write $X\oplus_1 Y$
and $X\oplus_\infty Y$ to denote, respectively, the $\ell_1$-sum
and the $\ell_\infty$-sum of $X$ and $Y$.

\section{Sufficient conditions for the Daugavet property}

For a better comprehension of the geometry underlying the Daugavet
property, we present the following characterization from
\cite[Lemma~2.1]{KSSW} and \cite[Corollary~2.3]{WerSur}. We shall
have occasion to use it throughout the paper.

\begin{lemma}\label{lemma-char}
The following assertions are equivalent:
\begin{equi}
\item $X$ has the Daugavet property.
\item For all $x\in S_X$, $x^*\in S_{X^*}$, and $\varepsilon>0$, there
exists some $y\in S_X$ such that $\re x^*(y)> 1- \varepsilon$ and
$\|x + y\|> 2-\varepsilon$.
\item For all $x\in S_X$, $x^*\in S_{X^*}$, and $\varepsilon>0$, there
exists some $y^*\in S_X$ such that $\re y^*(x)> 1- \varepsilon$
and $\|x^* + y^*\|> 2-\varepsilon$.
\item For all $x\in S_X$ and $\varepsilon>0$, $$B_X \subset
\ecc\bigl(\{y\in X\ : \ \|y\|\leqslant 1 + \varepsilon,\ \|x +
y\|> 2-\varepsilon\}\bigr).$$
\end{equi}
\end{lemma}

Observe that condition $(ii)$ in the above lemma implies that
every weak slice of the unit ball of a Banach space $X$ with the
Daugavet property has diameter $2$. Also, condition $(iii)$
implies that every $w^*$-slice of the unit ball of $X^*$ has
diameter $2$, thus the norm of the space is extremely rough.

The next result is a sufficient condition for a Banach space to
have the Daugavet property which will be crucial in the rest of
the paper. Recall that a closed subspace $Z$ of the dual of a
Banach space $X$ is called \emph{norming} whenever
$$\|x\|=\sup\{|z^*(x)|\ : \ z^*\in Z,\ \|z^*\|=1\} $$ for every
$x\in X$. This condition is clearly equivalent to $B_Z$ be
$w^*$-dense in $B_{X^*}$.

\begin{theorem}\label{th-sufDaugavet}
Let $X$ be a Banach space such that there are two norming
subspaces $Y$ and $Z$ of $X^*$ such that $X^{*}=Y\oplus_1 Z$.
Then, $X$ has the Daugavet property.
\end{theorem}

\begin{proof}
We fix $x_{0}\in S_{X}$, $f_0\in S_{X^*}$ and $\varepsilon>0$. We
write $f_0=y_0 + z_0$ such that $y_0\in Y$, $z_0\in Z$,
$\|f_0\|=\|y_0\| + \|z_0\|$, and $$U=\{x^*\in B_{X^*}\ : \ \re
x^*(x_0)>1 -\varepsilon\},$$ a $w^*$-open slice of $B_{X^*}$.
Since $B_{Z}$ is $w^*$-dense in $B_{X^*}$, we may find $z\in Z\cap
U$. Observe that, trivially, $\|z\|>1-\varepsilon$. Now, since
$B_Y$ is $w^*$-dense in $B_{X^*}$, we may find a net $(y_\lambda)$
in $B_Y$ which is $w^*$-convergent to $z$. Since $z\in U$, we may
suppose that $y_\lambda \in U$ for every $\lambda$. On the other
hand, since $(y_\lambda + y_0)\longrightarrow z + y_0$ and the
norm is $w^*$-lower semi-continuous, we have
$$\liminf \|y_\lambda + y_0\| \geqslant \|z +y_0\| = \|z\| + \|y_0\| >
1 + \|y_0\| -\varepsilon,$$ and we may find $\mu$ such that
$$\|y_\mu + y_0\| \geqslant 1 + \|y_0\| - \varepsilon/2.$$ To finish the
proof, we just observe that
\begin{align*}
\|f_0 + y_\mu\| & = \|(y_0 + y_\mu) + z_0\| \\ & = \|y_0 + y_\mu\|
+ \|z_0\| > 1 + \|y_0\| -\varepsilon + \|z_0\| = 2 -\varepsilon,
\end{align*}
and that $\re y_\mu(x_0)> 1-\varepsilon$ since $y_\mu\in U$, and
we use Lemma~\ref{lemma-char}.$iii$.
\end{proof}

Just remembering Goldstine and Krein-Milman Theorems, we obtain
the following useful particular case. Recall that a Banach space
$X$ is said to be $L$-embedded if $X^{**}=X\oplus_1 Z$ for some
closed subspace $Z$ of $X^{**}$.

\begin{corollary}\label{corollary1}
Let $X$ be a non-null $L$-embedded Banach space without extreme
points. Then, $X^*$ (and hence $X$) has the Daugavet property.
\end{corollary}

\begin{proof}
We have $X^{**}=X\oplus_1 Z$ for some subspace $Z$. On one hand,
since $B_{X}$ has no extreme points and
$\ext[X^{**}]=\ext[X]\,\bigcup \,\ext[Z]$, we have
$\ext[X^{**}]=\ext[Z]$ and the Krein-Milman Theorem gives us that
$B_Z$ is $w^*$-dense in $B_{X^{**}}$. On the other hand, Goldstine
Theorem gives us that $B_X$ is $w^*$-dense in $B_{X^{**}}$.
\end{proof}

It is worth mentioning that it is proved in \cite{Lop} that a
Banach space $X$ such that $X^{**}=X\oplus_1 Z$ with $B_Z$
$w^*$-dense in $B_{X^{**}}$ satisfies that every weak open subset
of $B_X$ has diameter two. Actually, the proof or
Theorem~\ref{th-sufDaugavet} has been inspired by the one given
there.

Let us finish the section by showing some immediate consequences
of the above result.

\begin{corollary}\label{cor-XYLembedded}
If $X$ is an $L$-embedded space with $\ext[X]=\emptyset$ and
$Y\subsetneq X$ is also an $L$-embedded space, then $(X/Y)^*$ (and
hence $X/Y$) has the Daugavet property.
\end{corollary}

\begin{proof}
$X/Y$ is a non-null $L$-embedded space
\cite[Corollary~IV.1.3]{HWW} and \cite[Propositions IV.1.12 and
IV.1.14]{HWW} gives us that $\ext[X/Y]=\emptyset$. Therefore,
Corollary~\ref{corollary1} applies.
\end{proof}

As a particular case of the above corollary we have the following
result.

\begin{corollary}
If $Y$ is an $L$-embedded subspace of $L_1\equiv L_1[0,1]$, then
$(L_1/Y)^*$ has the Daugavet property. In particular, $(L_1/Y)^*$
has the Daugavet property for every reflexive subspace $Y$ of
$L_1$ and so do $H^\infty$ and its predual $L_1/H_0^1$.
\end{corollary}

\begin{proof}
The space $L_1$ is an $L$-embedded space with
$\ext[L_1]=\emptyset$ and the space $H_0^1\subset L_1$ is also an
$L$-embedded space (see \cite[Example~IV.1.1]{HWW} for instance).
Then, the result follows immediately from
Corollary~\ref{cor-XYLembedded}.
\end{proof}

It is shown in \cite[Proposition~IV.2.11]{HWW} that $X/Y$ fails
the RNP when $X$ is an $L$-embedded space with $\ext[X]=\emptyset$
and $Y\subsetneq X$ is also an $L$-embedded space. On the other
hand, it is proved in \cite[Proposition~3.2]{KSSW} that $L_1/X$
has the Daugavet property whenever $X$ is a reflexive subspace of
$L_1$. The result for $H^\infty$ appeared in \cite{WerJFA} and
\cite{Woj}.

\section{$JB^*$-triples and preduals of $JBW^*$-triples}\label{sec-JB}

We recall that a \emph{complex $JB^*$-triple} is a complex Banach
space $X$ with  a continuous triple product $\{\cdots\}:\, X\times
X\times X\longrightarrow X$ which is linear and symmetric in  the
outer variables, and conjugate-linear in the middle variable, and
satisfies:
\begin{enumerate}
\item For all $x$ in $X$, the mapping $y\longmapsto \{xxy\}$ from
$X$ to $X$  is  a  hermitian operator on $X$ and has nonnegative
spectrum.
\item The \emph{main identity} $$\{ab\{xyz\}\}=\{\{abx\}yz\}-
\{x\{bay\}z\}+\{xy\{abz\}\} $$ holds for all $a,b,x,y,z$ in $X$.
\item $\Vert \{xxx\}\Vert =\Vert x\Vert ^{3}$ for every $x$ in
$X$.
\end{enumerate}
Concerning Condition (1) above, we also recall  that  a  bounded
linear operator $T$ on a complex Banach space $X$ is said to be
\emph{hermitian} if $\Vert \exp (irT)\Vert =1$ for every $r$ in
$\mathbb{R}$. By a \emph{complex $JBW^*$-triple} we mean a complex
$JB^*$-triple whose underlying Banach space is a dual space in
metric sense. It is known (see \cite{BaTi}) that every complex
$JBW^*$-triple has a unique predual up to isometric linear
isomorphisms and its triple product is separately $w^*$-continuous
in each variable.

Following \cite{IKR}, we define \emph{real $JB^*$-triples} as
norm-closed real subtriples of complex $JB^*$-triples. Here, by a
\emph{subtriple} we mean a subspace which is closed under triple
products of its elements. In particular, complex $JB^*$-triples
are real $JB^*$-triples. A \emph{triple ideal} of a real or
complex $JB^*$-triple $X$ is a subspace $M$ of $X$ such that
$\{XXM\}+ \{XMX\}\subseteq M$; if merely $\{M X M\}\subseteq M$,
then $M$ is called an \emph{inner ideal}.

\emph{Real $JBW^*$-triples} where first introduced as those real
$JB^*$-triples which are dual Banach spaces in such a way that the
triple product becomes separately $w^*$-continuous (see
\cite[Definition 4.1 and Theorem 4.4]{IKR}). Later, it has been
shown in \cite{MP} that the requirement of separate
$w^*$-continuity of the triple product is superabundant. We will
apply without notice that the bidual of every real or complex
$JB^*$-triple $X$ is a $JBW^*$-triple under a suitable triple
product which extends the one of $X$ (\cite{Dineen} for the
complex case and \cite{IKR} for the real case).

Examples of real $JB^*$-triples are the spaces $\mathcal{L}
(H,K)$, for arbitrary real, complex, or quaternionic Hilbert
spaces $H$ and $K$, under the triple product
$$\{xyz\}:=\frac{1}{2}(xy^*z+zy^*x).$$ The above examples become
particular cases of those arising by considering either the
so-called complex Cartan factors (regarded as real $JB^*$-triples)
or real forms of complex Cartan factors \cite{K10}. We recall that
\emph{real forms} of a complex Banach space $X$ are defined as the
real closed subspaces of $X$ of the form $X^\tau :=\{x\in X\ : \
\tau (x)=x\}$, for some conjugation (i.e., conjugate-linear
isometry of period two) on $X$. We note that, if $X$ is a complex
$JB^*$-triple, then every real form of $X$ is a real $JB^*$-triple
(since conjugations on $X$ preserve triple products \cite{K1}).
Conversely, if $X$ is a real $JB^*$-triple, there exists
\cite[Proposition~2.8]{IKR} a unique complex $JB^*$-triple
structure on the algebraic complexification $X\oplus\, i X$
(denoted $\widehat{X}$) and a conjugation $\tau$ on $X\oplus\, i
X$ such that $X=\widehat{X}^\tau$, i.e., every real $JB^*$-triple
is a real form of its complexification, which is a complex
$JB^*$-triple.

Let $X$ be a real or complex $JB^*$-triple. An element $u\in X$ is
said to be a \emph{tripotent} if $\{uuu\}=u$, and it said to be a
\emph{minimal tripotent} if $u\neq 0$ and $$\bigl\{x\in X\ : \
\{uxu\}=x\,\bigr\}=\mathbb{R} u.$$ In the complex setting, this is
equivalent to $u\neq 0$ and $\{uXu\}=\mathbb{C} u$.

If $x$ is a norm-one element of a real or complex $JB^*$-triple
$X$, then the set $D(X,x)=D(X^{**},x)\cap X^*$ is a proper closed
face of $B_{X^*}$, and therefore, by \cite[Lemma~2.1 and Theorem
3.7]{ER1}, there is a unique tripotent $u$ in $X^{**}$ such that
$D(X^{**},x)\cap X^*=D(X^{**},u)\cap X^*$. Such a tripotent $u$ is
called \emph{the support of $x$ in $X^{**}$}, and will be denoted
by $u(X^{**},x)$.

The complex case of the following result is stated in
\cite[Corollary~2.11]{BeRo}; the real case follows from results on
\cite{B-P} in an analogous way than the complex version. We
include the proof for the sake of completeness.

\begin{lemma}\label{lemma-FreMinTri}
Let $X$ be a real or complex $JB^*$-triple and let $x$ be in
$S_X$. Then, $X$ is Fr\'{e}chet-smooth at $x$ if and only if
$u(X^{**},x)$ lies in $X$ and it is a minimal tripotent of $X$.
\end{lemma}

\begin{proof}
Recall that the norm of a Banach space is Fr\'{e}chet-smooth at $x$ if
and only if it is smooth and strongly subdifferentiable at the
point (see \cite{FranPaya}). Now, the proof follows from the
following facts: the norm of $X$ is strongly subdifferentiable at
$x$ if and only if $u(X^{**},x)$ belongs to $X$
\cite[Corollary~2.5]{B-P}; $X$ is smooth at $x$ if and only if
$D(X^{**},x)\cap X^*=\{x^*\}$ for some extreme point $x^*$ of
$S_{X^*}$, and this is equivalent to the fact that $u(X^{**},x)$
is a minimal tripotent of $X^{**}$ \cite[Lemma~2.7 and
Corollary~2.1]{PeSta}; and, finally, a tripotent $u\in X$ is a
minimal tripotent of $X$ (if and) only if it is a minimal
tripotent of $X^{**}$ (by the $w^*$-density of $X$ in $X^{**}$ and
the separate $w^*$-continuity of the triple product of $X^{**}$).
\end{proof}

It is known \cite{BLPR} that the predual of every real or complex
$JBW^*$-triple is $L$-embedded. Therefore,
Corollary~\ref{corollary1} gives us that such a space has the
Daugavet property whenever its unit ball does not have any extreme
point. Actually, more can be proved:

\begin{theorem}\label{th-jbw}
Let $X$ be a real or complex $JBW^*$-triple and let $X_*$ be its
predual. Then, the following are equivalent:
\begin{equi}
\item $X$ has the Daugavet property.
\item $X_*$ has the Daugavet property.
\item Every relative weak-open subset of $B_{X_*}$ has diameter $2$.
\item $B_{X_*}$ has no strongly exposed points.
\item $B_{X_*}$ has no extreme points.
\end{equi}
\end{theorem}

\begin{proof}
$(i)\Rightarrow (ii)$ is clear. $(ii) \Rightarrow (iii)$ is
consequence of \cite[Lemma~3]{Shv}. $(iii) \Rightarrow (iv)$ is
clear.

$(iv)\Rightarrow (v)$. Of course, it is enough to show that every
extreme point of $B_{X_*}$ is actually an strongly exposed point.
Indeed, given $f\in\ext[X_*]$, \cite[Corollary~2.1]{PeSta} assures
the existence of a minimal tripotent $u$ of $X$ such that
$u(f)=1$, and $u$ is a point of Fr\'{e}chet-smoothness of the norm of
$X$ by Lemma~\ref{lemma-FreMinTri}. Therefore, there is a point of
Fr\'{e}chet-smoothness, $u$, in $D(X_*,f)$ and, as we commented in the
introduction, this implies that $f$ is strongly exposed by $u$
(see \cite[Corollary~I.1.5]{DGZ}, for instance).

$(v)\Rightarrow (i)$. $X_*$ is an L-embedded by
\cite[Proposition~2.2]{BLPR} and $B_{X_*}$ has no extreme points,
so Corollary~\ref{corollary1} applies.
\end{proof}

As an straightforward consequence of the above theorem we obtain
the following result, which states the ``extreme'' behaviour of
the diameters of the weak-open subset of the unit ball of the
predual of a $JBW^*$-triple.

\begin{corollary}\label{cor-jbweither}
Let $Y$ be the predual of some real or complex $JBW^*$-triple.
Then, either every weak-open subset of $B_Y$ has diameter $2$ or
$B_Y$ has slices of arbitrary small diameter.
\end{corollary}

By Corollary~2.1 of \cite{PeSta}, a real or complex $JBW^*$-triple
has minimal tripotents if and only if the unit ball of its predual
has extreme points. Therefore, the following result follows
immediately from Theorem~\ref{th-jbw}.

\begin{corollary}\label{corollary-jbw}
Let $X$ be a real or complex $JBW^*$-triple. Then, $X$ has the
Daugavet property if and only if it does not have any minimal
tripotents.
\end{corollary}

The complex case of the above corollary and the equivalence
$(i)\Leftrightarrow (ii)$ of Theorem~\ref{th-jbw} appear in
\cite[Theorem~4.7]{MaOi}.

As a consequence of Theorem~\ref{th-jbw} we obtain:

\begin{corollary}\label{cor-nordualnorbidual}
Neither the dual of a real or complex $JB^*$-triple nor a real or
complex $JB^*$-triple which is the bidual of some space, has the
Daugavet property.
\end{corollary}

\begin{proof}
On one hand, the dual $X^*$ of a $JB^*$-triple $X$ is also the
predual of the $JBW^*$-triple $X^{**}$ and, as every dual space,
$B_{X^*}$ has extreme points. On the other hand, if $Y=Z^{**}$ is
a $JB^*$-triple, then it is actually a $JBW^*$-triple whose
predual $Y_*=Z^*$ has extreme points in its unit ball.
\end{proof}

\begin{remark}
It is worth mentioning that, for an arbitrary Banach space $Z$,
the absence of extreme points in $B_Z$ or the fact that all
weak-open subsets of $B_Z$ have diameter two, does not necessarily
imply that $Z$ has the Daugavet property. For instance, $c_0$
satisfies both assumptions (see \cite[Lemma~2.2]{BeLoRo} for
instance),  but it does not have the Daugavet property.

On the other hand, the assertions $(iii)$, $(iv)$, and $(v)$ of
Theorem~\ref{th-jbw} are not equivalent for general Banach spaces.
On one hand, {\slshape there exists a Banach space $Z$ whose unit
ball has slices of arbitrary small diameter, but it does not have
any extreme point (so, it does not have any strongly exposed
point) \cite[Proposition~1]{Edel}.\ } On the other hand, {\slshape
every slice of the unit ball of $\ell_\infty$ has diameter $2$
(and so, it does not have any strongly exposed point), but it is
plenty of extreme points (it is a dual space).}
\end{remark}

If $X$ is a real or complex $JBW^*$-triple, it is well known that
$X_*=A\oplus_1 N$, where $A$ is the closed linear span of the
extreme points of $B_{X_*}$, and the unit ball of $N$ has no
extreme points (see \cite{FR85} for the complex case and
\cite{PeSta} for the real case). Therefore,
$X=\mathcal{A}\oplus_\infty \mathcal{N}$, where
$\mathcal{A}=N^\perp\equiv A^*$ is an \emph{atomic} $JBW^*$-triple
(i.e.\ it is the weak*-closed span of its minimal tripotents) and
$\mathcal{N}=A^\perp\equiv N^*$ is a $JBW^*$-triple without
minimal tripotents. With this in mind, the following result is a
consequence of Theorem~\ref{th-jbw} and a characterization of the
RNP in preduals of $JBW^*$-triples given in \cite{BaGo}.

\begin{corollary}\label{cor-jbwRNP-Daugavet}
Let $X$ be a real or complex $JBW^*$-triple. Then, in the natural
decomposition  $X_*=A\oplus_1 N$, $A$ has the RNP and $N$ has the
Daugavet property. Therefore, in the decomposition
$X=\mathcal{A}\oplus_\infty \mathcal{N}$, $\mathcal{A}$ is a
$w^*$-Asplund space (i.e., the dual of a space having the RNP) and
$\mathcal{N}$ has the Daugavet property.
\end{corollary}

\begin{proof}
In the complex case, since $A$ is the predual of the atomic
$JBW^*$-triple $\mathcal{A}$, it has the RNP by
\cite[Theorem~1]{BaGo} and, therefore, $\mathcal{A}$ is a
$w^*$-Asplund space. In the real case, we consider
$\widehat{\mathcal{A}}$, the complexification of $\mathcal{A}$. On
one hand, $\widehat{\mathcal{A}}$ is a $w^*$-Asplund space by the
above. On the other hand, $A\equiv \mathcal{A}_*$ is a (real)
subspace of $\bigl(\widehat{\mathcal{A}}\,\bigr)_*$, and the RNP
passes to subspaces.

Since $N^*=\mathcal{N}$ is a $JBW^*$-triple without minimal
tripotents, Corollary~\ref{corollary-jbw} gives us that
$\mathcal{N}$, and hence its predual $N$, have the Daugavet
property.
\end{proof}

Our next aim is to prove a characterization of the Daugavet
property for general $JB^*$-triples. We first prove that the
algebraic characterization given in Corollary~\ref{corollary-jbw}
for $JBW^*$-triples is also valid in the general case, and then we
will deduce more characterizations in terms of the geometry of the
norm of the triple.

We need a result about real or complex $JB^*$-triples which can be
of independent interest. Previously, we have to recall some known
facts about $JB^*$-triples.

If $X$ is a real or complex $JB^*$-triple, $X^{**}$ is a
$JBW^*$-triple. Therefore, we can decompose $X^{*}=(X^{**})_*$
into its atomic and not atomic parts, as we have commented above,
i.e., $X^*=A\oplus_1 N$ where $A$ is the closed linear span of the
extreme points of $B_{X^*}$, and the unit ball of $N$ has no
extreme points. Then, $X^{**}=\mathcal{A} \oplus_\infty
\mathcal{N}$, where $\mathcal{A}=N^\perp\equiv A^*$ is an atomic
$JBW^*$-triple, and $\mathcal{N}=A^\perp\equiv N^*$ is a
$JBW^*$-triple without minimal tripotents. Let us call
$\pi_\mathcal{A}$ (resp.\ $\pi_\mathcal{N}$) the projection from
$X^{**}$ to $\mathcal{A}$ with kernel $\mathcal{N}$ (resp.\ to
$\mathcal{N}$ with kernel $\mathcal{A}$), and let
$J_X:X\longrightarrow X^{**}$ be the natural inclusion. It is well
known that $\pi_\mathcal{A}\circ J_X:X \longrightarrow
\mathcal{A}$ is an isometric embedding (Gelfand-Naimark Theorem
\cite{FR86}). The next result gives the same for
$\pi_\mathcal{N}\circ J_X$, provided $X$ has no minimal
tripotents.

\begin{theorem}\label{th-GN2}
Let $X$ be a real or complex $JB^*$-triple without minimal
tripotents. Then, the mapping $\pi_\mathcal{N}\circ J_X:X
\longrightarrow \mathcal{N}$ is an isometric embedding. Therefore,
$N$ is a norming subspace of $X^*$.
\end{theorem}

\begin{proof}
We start by proving the result in the complex case. Let $X$ be a
complex $JB^*$-triple and let us consider $Y=X\cap \mathcal{A}$,
which is clearly a closed ideal of $X$. On one hand, $Y$ has no
minimal tripotents (indeed, if $0\neq u\in Y$ is a minimal
tripotent of $Y$, then $\{u Y u\}=\mathbb{C} u$; since $Y$ is a
triple ideal (and hence an inner ideal), we have $\{u X u\}\subset
Y$, so we obtain $\{u X u\}=\mathbb{C} u$ and $u$ is a minimal
tripotent of $X$, which is impossible). On the other hand, by
\cite[Proposition~3.7]{BC91} $Y^*$ has the RNP (i.e.\ $Y$ is an
Asplund space) and, if $Y\neq 0$, the norm of $Y$ has points of
Fr\'{e}chet-smoothness. But the existence of points of
Fr\'{e}chet-smoothness in $Y$ implies the existence of minimal
tripotents in $Y$ (Lemma~\ref{lemma-FreMinTri}), a contradiction.
We deduce that $Y$ is null and, therefore, $\pi_\mathcal{N}\circ
J_X$ is injective. Being a triple-homomorphism, it is routine
(using axiom (3)) to show that it is an isometric embedding as
desired (actually, in the complex case, the converse result is
also true, see \cite{K1}). Since $\mathcal{N}=A^\perp\equiv N^*$,
it is clear that $N$ is norming.

The proof for the real case is very similar. If $X$ is a real
$JB^*$-triple, we will show that $Y=X\cap \mathcal{A}$ has no
minimal tripotents and that it is an Asplund space, and then the
rest of the above proof works. First, if $0\neq u \in Y$ is a
minimal tripotent, then $\{y\in Y\ : \ \{u y u\}= y\}=\mathbb{R}
u$; since $Y$ is a inner ideal, $\{u X u\}\subseteq Y$, so if
$x\in X$ is such that $\{u x u\}=x$, we obtain that $x\in Y$,
which implies $x\in \mathbb{R} u$, i.e., $u$ is a minimal
tripotent of $X$, a contradiction. Second, we consider the
complexification $\widehat{Y}$ of $Y$, and we observe that
$\widehat{Y}=\widehat{\mathcal{A}\,}\,\cap \widehat{X\,}$, where
$\widehat{X^{**}\,}=\widehat{\mathcal{A}\,}\,\oplus_\infty
\widehat{\mathcal{N}\,}$ is the decomposition into the atomic and
non-atomic part \cite[Theorem~3.6]{PeSta}. Therefore,
$\widehat{Y}$ is an Asplund space \cite[Proposition~3.7]{BC91} and
so does its real subspace $Y$.
\end{proof}

As a consequence of the above result and
Theorem~\ref{th-sufDaugavet}, we obtain that $JB^*$-triples
without minimal tripotents have the Daugavet property. The complex
case of this result appear in \cite[Theorem~4.7]{MaOi} with a
different proof.

\begin{prop}\label{prop-JB}
Let $X$ be a real or complex $JB^*$-triple. Then, $X$ has the
Daugavet property if and only if it has no minimal tripotents.
\end{prop}

\begin{proof}
Suppose $X$ has no minimal tripotents and write $X^*=A\oplus_1 N$.
On one hand, since $\ext \subseteq B_A$, the Krein-Milman Theorem
gives us that $A$ is a norming subspace of $X^*$. On the other
hand, if $X$ has no minimal tripotents, Theorem~\ref{th-GN2} gives
us that $N$ is also norming. Now, Theorem~\ref{th-sufDaugavet}
gives us that $X$ has the Daugavet property. Conversely, if $X$
has a minimal tripotents, then it has a point of
Fr\'{e}chet-smoothness by Lemma~\ref{lemma-FreMinTri}; but the norm of
a Banach space with the Daugavet property is extremely rough (use
Lemma~\ref{lemma-char}.$iii$), a contradiction.
\end{proof}

Actually, we can state a characterization of the Daugavet property
for $JB^*$-triples in terms of the geometry of the norm of the
triple.

\begin{theorem}\label{th-JB}
Let $X$ be a real or complex $JB^*$-triple. Then, the following
are equivalent:
\begin{equi}
\item $X$ has the Daugavet property.
\item The norm of $X$ is extremely rough.
\item The norm of $X$ is not Fr\'{e}chet-smooth at any point.
\end{equi}
\end{theorem}

\begin{proof}
$(i)\Rightarrow (ii)$. As we commented in the introduction, the
norm of $X$ is extremely rough if and only if every $w^*$-slice of
$B_{X^*}$ has diameter 2, and the latest fact is consequence of
Lemma~\ref{lemma-char}.$iii$.

$(ii) \Rightarrow (iii)$ is clear.

$(iii)\Rightarrow (i)$. By Lemma~\ref{lemma-FreMinTri}, the norm
of $X$ is Fr\'{e}chet-smooth at the minimal tripotents, so we deduce
that $X$ has no minimal tripotents and Proposition~\ref{prop-JB}
applies.
\end{proof}

\begin{remark}
It is worth mentioning that the above geometric characterizations
are not valid for arbitrary Banach spaces. For instance, {\slshape
the norm of $\ell_1$ is extremely rough (and so $\ell_1$ has no
points of Fr\'{e}chet-smoothness) but $\ell_1$ does not have the
Daugavet property.}

Also, the implication $(iii)\Rightarrow (ii)$ of the above theorem
is not valid in general. Indeed, {\slshape there exists a Banach
space whose norm does not have any point of Fr\'{e}chet
differentiability but it is not rough} (see \cite[Remark~4,
pp.~341]{JoZi}).
\end{remark}

To finish the section, let us comment some results from
\cite{BLPR} which are related to our development.

\begin{remark}
Let us consider the following conditions for a Banach space $X$:
\begin{enumerate}
\item[(a)] every relative weak-open subset of $B_X$ has diameter
$2$,
\item[(b)] the norm of $X$ is extremely rough.
\end{enumerate}

It is proved in \cite[Theorem~2.3]{BLPR} that condition (a) is
satisfied when $X$ is a non-reflexive real or complex
$JB^*$-triple, while our Theorem~\ref{th-jbw} says that condition
(a) characterizes the Daugavet property in the class of preduals
of real or complex $JBW^*$-triples.

With respect to condition (b), it is shown in
\cite[Corollary~2.5]{BLPR} that the predual of every non-reflexive
real or complex $JBW^*$-triple satisfies it, while condition (b)
characterizes the Daugavet property for real or complex
$JB^*$-triples (Theorem~\ref{th-JB}).

Since a reflexive Banach space never satisfies neither (a) nor
(b), the above paragraphs contains the answer to every question
about this conditions in the setting of real or complex
$JB^*$-triples and their isometric preduals.
\end{remark}

\section{$C^*$-algebras and von Neumann preduals}

Despite \emph{real $C^*$-algebras} can be defined by different
systems of intrinsic axioms (see \cite{IsiRod} for a summary), we
prefer to introduce them as the norm-closed self-adjoint real
subalgebras of complex $C^*$-algebras. Since complex
$C^*$-algebras are complex $JB^*$-triples under the triple product
$$\{xyz\}:=\frac{1}{2}(xy^*z+zy^*x),$$ certainly real $C^*$-algebras
are real $JB^*$-triples. The concept of a \emph{real
$W^*$-algebra} (\emph{real von Neumann algebra}) was first defined
as a real $C^*$-algebra $A$ having a complete predual $A_*$ such
that the product of $A$ is separately $w^*$-continuous, but the
latest condition was shown to be redundant in \cite{IsiRod}. Real
$W^*$-algebras are real $JBW^*$-triples.

Therefore, the geometric characterizations given in Theorems
\ref{th-jbw} and \ref{th-JB} can be stated for real or complex
$C^*$-algebras and preduals of $W^*$-algebras. The next results
summarize those theorems and also Corollaries \ref{cor-jbweither}
and \ref{cor-nordualnorbidual} in terms of $C^*$-algebras.

\begin{corollary}
Let $X$ be a real or complex $C^*$-algebra. Then, the following
are equivalent:
\begin{equi}
\item $X$ has the Daugavet property.
\item The norm of $X$ is extremely rough.
\item The norm of $X$ is not Fr\'{e}chet-smooth at any point.
\end{equi}
\end{corollary}

\begin{corollary}
Let $X$ be a real or complex $W^*$-algebra and let $X_*$ be its
predual. Then, the following are equivalent:
\begin{equi}
\item $X$ has the Daugavet property.
\item $X_*$ has the Daugavet property.
\item Every weak-open subset of $B_{X_*}$ has diameter $2$.
\item $B_{X_*}$ has no strongly exposed points.
\item $B_{X_*}$ has no extreme points.
\end{equi}
\end{corollary}

\begin{corollary}$ $
\begin{enumerate}
\item[(a)] Let $X$ be the predual of some real or complex
$W^*$-algebra. Then, either every weak-open subset of $B_X$ has
diameter $2$ or $B_X$ has slices of arbitrary small diameter.
\item[(b)] Neither the dual of a real or complex $C^*$-algebra nor a real
or complex $C^*$-algebra which is the bidual of some space, has
the Daugavet property.
\end{enumerate}
\end{corollary}

The algebraic characterization of the Daugavet property for
$JB^*$-triples (Proposition~\ref{prop-JB}) is of course valid for
$C^*$-algebras, but it could be more convenient to write it in
terms of atomic projections. Let us gives the definitions and
results.

If $X$ is a real or complex $C^*$-algebra, then $u \in X$ is a
tripotent if and only if it is a \emph{partial isometry}, i.e.,
$u$ satisfies that $uu^*u=u$. Recall that a \emph{projection} in a
$C^*$-algebra is an element $p\in X$ such that $p^*=p$ and
$p^2=p$. It is clear that projections are partial isometries (and
so tripotents), but there are partial isometries which are not
projections. A projection $p$ in $X$ is said to be \emph{atomic}
if $p\neq 0$ and $$\left\{x\in X\ : \ px^*p=x\right\}=\mathbb{R}
p,$$ i.e., $p$ is minimal seen as a tripotent. Therefore, in the
complex case this is equivalent to $p\neq 0$ and $p X p=\mathbb{C}
p$. The $C^*$-algebra $X$ is said to be \emph{non-atomic} if it
does not have any atomic projection.

If $X$ has atomic projections, then it clearly has minimal
tripotents. Conversely, if $X$ has a minimal tripotent, say $u$,
then the projection $d=u^* u$ (called the domain projection
associated to $u$) is atomic. Indeed, we take $x\in X$ such that
$dx^*d=x$. Then,
$$u(ux)^*u=(uu^*u) x^* u^* u = u (u^*u x^* u^* u) =u (dx^*d) = u x$$
so, since $u$ is minimal, $ux=\lambda u$ for some $\lambda\in
\mathbb{R}$. Then,
\begin{align*}
\lambda d =u^*(\lambda u) &=u^*(ux)=u^*(u(dx^*d))=
\\ & =u^*\bigl((uu^*u)x^*u^*u\bigr)=u^*ux^*u^*u=dx^*d=x.
\end{align*}
We have shown that {\slshape a real or complex $C^*$-algebra has
no minimal tripotents if and only if it is non-atomic.\ } So, for
$C^*$-algebras, Proposition~\ref{prop-JB} can be written in terms
of atomic projections.

\begin{corollary}$ $
\begin{enumerate}
\item[(a)] A real or complex $C^*$-algebra has the Daugavet property
if and only if it is non-atomic.
\item[(b)] The predual of a real or complex $W^*$-algebra has the Daugavet
property if and only if the algebra is non-atomic.
\end{enumerate}
\end{corollary}

The complex case of the above result appear in
\cite[Theorem~2.1]{Oik}.

As a $JBW^*$-triple, every real or complex $W^*$-algebra $X$
admits a natural decomposition into the atomic and non-atomic
parts which is originated by the natural decomposition of the
predual $X_*$. I.e., $X_*=A\oplus_1 N$, where the unit ball of $N$
does not have any extreme point, and $B_A$ is the closed convex
hull of the extreme points of $B_{X_*}$. Thus,
$X=\mathcal{A}\oplus_\infty \mathcal{N}$, where the subtriple
$\mathcal{A}=N^\perp\equiv A^*$ is norm-generated by the minimal
tripotents of $X$, and the subtriple $\mathcal{N}=A^\perp\equiv
N^*$ has no minimal tripotents. Moreover, $\mathcal{A}$ and
$\mathcal{N}$ are $w^*$-closed subalgebras of $X$, the first one
is generated by its atomic projections and the second one has no
atomic projections.

The next results put Corollary~\ref{cor-jbwRNP-Daugavet} and
Theorem~\ref{th-GN2} in terms of $C^*$-algebras.

\begin{corollary}
Let $X$ be a real or complex $W^*$-algebra. Then, in the natural
decomposition  $X_*=A\oplus_1 N$, $A$ has the RNP and $N$ has the
Daugavet property. Therefore, in the decomposition
$X=\mathcal{A}\oplus_\infty \mathcal{N}$, $\mathcal{A}$ is a
$w^*$-Asplund space (i.e., the dual of a space having the RNP) and
$\mathcal{N}$ has the Daugavet property.
\end{corollary}

\begin{corollary}
Let $X$ be a real or complex $C^*$-algebra without atomic
projections, and let $X^{**}=\mathcal{A}\oplus_\infty \mathcal{N}$
the natural decomposition of its bidual into atomic and non-atomic
parts. Then, the decomposition of every $x\in X$ as $x=a^{**} +
n^{**}$, with $a^{**}\in \mathcal{A}$, $n^{**}\in\mathcal{N}$
satisfies $\|x\|=\|a^{**}\|=\|n^{**}\|$.
\end{corollary}

\section{The uniform Daugavet property}
Following \cite{BKSW}, a Banach space $X$ is said to have the
\emph{uniform Daugavet property} if $$D_X(\varepsilon):=\inf
\{n\in \mathbb{N}\ : \ \text{conv}_n(l^+(x,\varepsilon))\supset
S_X\ \forall x\in S_X\}$$ is finite for every $\varepsilon>0$,
where
$$l^+(x,\varepsilon):=\bigl\{y\in X\ : \ \|y\|\leqslant 1 + \varepsilon,\ \|x + y\|>
2-\varepsilon\bigr\}$$ and $\text{conv}_n(A)$ is the set of all
convex combination of all $n$-point collections of elements of
$A$. By \cite[Remark~6.3]{BKSW}, $X$ has the uniform Daugavet
property if and only if
$$\lim_{n\to\infty}\text{Daug}_{\,n}(X,\varepsilon) = 0$$ for every
$\varepsilon>0$, where
$$\text{Daug}_{\,n}(X,\varepsilon):=\  \sup_{x,y\in S_X} \
\text{dist}\bigl(y,\text{conv}_n(l^+(x,\varepsilon))\bigr).$$
Since (Lemma~\ref{lemma-char}) $X$ has the Daugavet property if
and only if $$B_X \subset \ecc\bigl(\{y\in X\ : \ \|y\|\leqslant 1
+ \varepsilon,\ \|x + y\|> 2-\varepsilon\}\bigr)$$ for every $x\in
S_X$ and every $\varepsilon>0$, the uniform Daugavet property
implies the Daugavet property, and it can be view as a
quantitative approach to it.

Examples of spaces satisfying the uniform Daugavet property are
$L_1[0,1]$ and $C(K)$ for every perfect compact space $K$ \cite[\S
6]{BKSW}. On the other hand, in \cite{KaWe} it is shown an example
of a Banach space with the Daugavet property which does not
satisfy the uniform Daugavet property.

The uniform Daugavet property was introduced in \cite{BKSW} to
study when the Daugavet property passes from a Banach space to its
so-called ultrapowers.

Let us recall here the notion of (Banach) ultraproducts \cite{He}.
Let $\mathcal{U}$ be a free ultrafilter on a nonempty set $I$, and
let $\{ X_i \}_{i\in I}$ be a family of Banach spaces. We can
consider the $\ell_\infty$-sum of the family, $\left[\oplus_{i\in
I} X_i\right]_{\ell_\infty}$, together with its closed subspace
$$N_{\mathcal U} := \left\{ \{x_i\}_{i\in I} \in \left[\oplus_{i\in
I} X_i\right]_{\ell_\infty} \ :\  \lim_{\mathcal U} \Vert x_i
\Vert=0 \right\}.$$ The quotient space $\left[\oplus_{i\in I}
X_i\right]_{\ell_\infty}/ N_{\mathcal U}$ is called the
\emph{ultraproduct} of the family $\{ X_i \}_{i\in I}$ relative to
the ultrafilter ${\mathcal U}$, and is denoted by $(X_i)_{\mathcal
U}$. Let $(x_i)$ stand for the element of $(X_i)_{\mathcal U}$
containing a given family $\{ x_i \} \in \left[\oplus_{i\in I}
X_i\right]_{\ell_\infty}$. It is easy to check that $\Vert (x_i)
\Vert = \lim_{\mathcal U} \Vert x_i \Vert$. Moreover, the
ultraproduct $(X_i^*)_{\mathcal U}$ can be seen as a subspace of
$\bigl[(X_i)_{\mathcal U}\bigr]^*$ by identifying each element
$(f_i)\in (X_i^*)_{\mathcal U}$ with the (well-defined) functional
on $(X_i)_{\mathcal U}$ given by $$(x_i)\longmapsto \lim _\mathcal
U (f_i(x_i)) \qquad \bigl((x_i)\in (X_i)_\mathcal{U}\bigr).$$ If
$\{Y_i\}_{i\in I}$ is another family of Banach spaces and for each
$i\in I$ we take an operator $T_i\in L(X_i,Y_i)$ with $\sup_{i\in
I}\|T_i\|<\infty$, we can define the \emph{utraproduct of the
family of operators} $\{T_i\}_{i\in I}$ with respect to the
ultrafilter $\mathcal{U}$, denoted $(T_i)$, as
$$(x_i)\longmapsto (T_i x_i)\qquad
\bigl((x_i)\in (X_i)_\mathcal{U}\bigr).$$ This is a well defined
operator from $(X_i)_\mathcal{U}$ to $(Y_i)_\mathcal{U}$ with
$$\left\|(T_i)\right\| = \lim_\mathcal{U}\|T_i\|.$$ If
all the $X_i$ are equal to some Banach space $X$, the ultraproduct
of the family is called the $\mathcal{U}$-\emph{ultrapower} of $X$
and it is usually denoted by $X_\mathcal{U}$. For $T\in L(X)$, by
$(T)$ we denote the ultraproduct of the family $\{T_i\}_{i\in I}$
where $T_i=T$ for every $i\in I$.

In \cite[Corollary~6.5]{BKSW}, it is proved that a Banach space
$X$ has the uniform Daugavet property if and only if every
ultrapower $X_\mathcal{U}$, $\mathcal{U}$ a free ultrafilter on
$\mathbb{N}$, has the Daugavet property, in which case
$X_\mathcal{U}$ even has the uniform Daugavet property. Let us
comment that it is routine to prove that a Banach $X$ has the
(usual) Daugavet property whenever $X_\mathcal{U}$ does,
$\mathcal{U}$ a free ultrafilter on an arbitrary set $I$ (we can
use Lemma~\ref{lemma-char}.$ii$ or, alternatively, we can prove
directly that every rank-one operator $T\in L(X)$ satisfies
\eqref{DE} since its ultrapower $(T)\in L(X_\mathcal{U})$, which
is also a rank-one operator on $X_\mathcal{U}$, does). On the
other hand, as we have said before, there is a Banach space with
the Daugavet property which does not have the uniform Daugavet
property \cite{KaWe}, thus the Daugavet property does not always
pass to ultrapowers.

Our aim in this section is to prove that the Daugavet property and
its uniform version are equivalent for real or complex
$JB^*$-triples and their isometric preduals. As we said before,
this is true for $C(K)$ spaces and for $L_1[0,1]$. These facts
were proved in \cite[\S 6]{BKSW}, where explicit estimations for
$D_{C(K)}(\varepsilon)$ and $D_{L_1[0,1]}(\varepsilon)$ were done.
Our approach is different: we will use Theorems \ref{th-jbw} and
\ref{th-JB} to show that, for $JB^*$-triples and their isometric
preduals, the Daugavet property passes to arbitrary ultrapowers.

Since an ultrapower of a $JB^*$-triple is again a $JB^*$-triple
(see \cite{Dineen}), the result for this class follows immediately
from Theorem~\ref{th-JB} and the following lemma, which can be of
independent interest.

\begin{lemma}\label{lemmaroughtoultraproducts}
Let $\{ X_i \}_{i\in I}$ be a family of Banach spaces,
$\mathcal{U}$ a free ultrafilter of a set $I$, and $\delta>0$. If
the norm of each $X_i$ is $\delta$-rough, then so does the norm of
$(X_i)_{\mathcal U}$.
\end{lemma}

\begin{proof}
Given a norm-one element $x=(x_i)\in (X_i)_\mathcal{U}$ and a
positive number $\alpha<1$, we have to show that the slice
$S\bigl(B_{\left[(X_i)_\mathcal{U}\right]^*},(x_i),\alpha\bigr)$
of the unit ball of $\bigl[(X_i)_\mathcal{U}\bigr]^*$ has diameter
greater than $\delta$. Indeed, we can suppose that $\|x_i\|=1$ for
every $i\in I$ and, since the norm of each $X_i$ is
$\delta$-rough, given a family $\{\varepsilon_i\}$ of positive
number with $\lim_\mathcal{U}\varepsilon_i=0$, we can find
$f_i,g_i\in S_{X_i^*}$ such that
$$\|f_i-g_i\| > \delta -\varepsilon_i \qquad \text{and} \qquad \re f_i(x_i)
> 1 - \alpha,\quad \re g_i(x_i)> 1 - \alpha.$$ Now, we consider the
elements $f=(f_i)$ and $g=(g_i)$ of the unit ball of
$(X_i^*)_\mathcal{U}\subseteq \bigl[(X_i)_\mathcal{U}\bigr]^*$,
and we observe that, on one hand, $$\bigl\|(f_i)-(g_i)\bigr\|
=\lim_\mathcal{U}\|f_i-g_i\|\geqslant \delta$$ and, on the other
hand,
\begin{equation*}
\re f(x) =\lim_{\mathcal{U}}f_i(x_i)>1 - \alpha, \ \qquad \re g(x)
=\lim_{\mathcal{U}}g_i(x_i)>1 - \alpha.\qedhere
\end{equation*}
\end{proof}

By using the above lemma and Theorem~\ref{th-JB}, we have that
$X_\mathcal{U}$ has the Daugavet property whenever the
$JB^*$-triple $X$ does. But, as we already mentioned, the converse
result is true in general.

\begin{theorem}\label{th-ultraproductJB}
Let $X$ be a real or complex $JB^*$-triple and $\mathcal{U}$ a
free ultrafilter on a set $I$. Then, $X$ has the Daugavet property
if and only if $X_\mathcal{U}$ does. Therefore, the Daugavet
property and the uniform Daugavet property are equivalent for
$JB^*$-triples.
\end{theorem}

As a consequence of the above theorem and
Proposition~\ref{prop-JB}, we obtain the following result about
$JB^*$-triples.

\begin{corollary}
Let $X$ be a real or complex $JB^*$-triple and $\mathcal{U}$ a
free ultrafilter on a set $I$. Then, $X_\mathcal{U}$ has a minimal
tripotent if and only if $X$ does.
\end{corollary}

\begin{remark}
It is also true that every ultraproduct of $JB^*$-triples is a
$JB^*$-triple (see \cite{Dineen}). Then, by using
Theorem~\ref{th-JB} and Lemma~\ref{lemmaroughtoultraproducts}, we
also obtain that {\slshape the ultraproduct of a family of
$JB^*$-triples with the Daugavet property also has the Daugavet
property.\ } In other words (Proposition~\ref{prop-JB}), {\slshape
the ultraproduct of a family of $JB^*$-triples without minimal
tripotents also has no minimal tripotent.}
\end{remark}

The second part of the present section is devoted to preduals of
$JBW^*$-triples.

Even though the ultrapower of the dual of a Banach space is not,
in general, the dual of the ultrapower of the space (see \cite[\S
7]{He}), it can be proved that the ultrapower of a predual of a
$JBW^*$-triple is again the predual of some $JBW^*$-triple. In the
complex case, the proof is easy to state: the dual of the
ultrapower $X_\mathcal{U}$ of a Banach space $X$ is
$1$-complemented in another ultrapower $(X^*)_\mathcal{M}$ of
$X^*$ \cite{He}, and the contractive projection Theorem applies.

Since we have not find any reference to the above result in the
literature, we give a detailed proof. Actually, a more general
result can be state.

\begin{prop}\label{prop-preUltraproducts}
Let $\{X_i\}_{i\in I}$ a family of Banach spaces such that each
$X_i^*$ is a (real or complex) $JBW^*$-triple, and let
$\mathcal{U}$ be a free ultrafilter on $I$. Then,
$(X_i)_\mathcal{U}$ is the predual of some (real or complex)
$JBW^*$-triple.
\end{prop}

\begin{proof}
We start with the complex case. By \cite[Corollary~7.6]{He}, there
is another free ultrafilter $\mathcal{B}$ on an index set $I'$,
such that $\bigl[(X_i)_\mathcal{U}\bigr]^*$ is isometric to a
$1$-complemented subspace of $((X_i^*)_\mathcal{U})_\mathcal{B}$,
which is a $JB^*$-triple. But $1$-complemented subspaces of
complex $JB^*$-triples are $JB^*$-triples (see \cite{Kau84}).

If each $X_i^*$ is a real $JBW^*$-triple, then there is a
conjugation $\tau_i$ on each $X_i$ such that
$\bigl(\widehat{X_i}\bigr)^*$ is a complex $JBW^*$-triple and
$X_i=\widehat{X_i}^{\tau_i}$ \cite{IKR}. On one hand,
$\bigl[\bigl(\widehat{X_i}\bigr)_\mathcal{U}\bigr]^*$ is a
$JBW^*$-triple by the complex case. On the other hand, we consider
$\tau=(\tau_i)$, the ultraproduct of the family of the
conjugations $\tau_i$, and we observe that $\tau$ is a conjugation
(routine) and that $\bigl[(\widehat{X_i})_\mathcal{U}\bigr]^\tau
\equiv (X_i)_\mathcal{U}.$ Indeed, $(\widehat{x_i})\in
\bigl[(\widehat{X_i})_\mathcal{U}\bigr]^\tau$ if and only if
$\lim_\mathcal{U}\|\tau_i(\widehat{x_i})-\widehat{x_i}\|=0$. Thus,
the image of the natural inclusion of $(X_i)_\mathcal{U}$ into
$(\widehat{X_i})_\mathcal{U}$ falls into
$\bigl[(\widehat{X_i})_\mathcal{U}\bigr]^\tau$, and it is onto
since, for every $(\widehat{x_i})\in
\bigl[(\widehat{X_i})_\mathcal{U}\bigr]^\tau$, we have
$(\widehat{x_i})= \bigl(\tau_i(\widehat{x_i})\bigr)\in
\bigl(\widehat{X_i}^{\tau_i}\bigr)_\mathcal{U}\equiv
(X_i)_\mathcal{U}$. Now, the dual of $(X_i)_\mathcal{U}\equiv
\bigl[(\widehat{X_i})_\mathcal{U}\bigr]^\tau$ is a real form
(using $\tau^*$, which is also a conjugation) of
$\bigl[(\widehat{X_i})_\mathcal{U}\bigr]^*$, and hence it is a
real $JBW^*$-triple.
\end{proof}

With this in mind, the equivalence of the Daugavet property and
its uniform version for preduals of $JBW^*$-triples is a
consequence of Theorem~\ref{th-jbw}.

\begin{theorem}\label{th-DaugtoultrapowerpredulaJBW}
Let $X$ be a real or complex $JBW^*$-triple and $\mathcal{U}$ a
free ultrafilter on a set $I$. Then, $X_*$ has the Daugavet
property if and only if $(X_*)_\mathcal{U}$ does. Therefore, the
Daugavet property and the uniform Daugavet property are equivalent
for preduals of $JBW^*$-triples.
\end{theorem}

In the proof we will use the following easy fact: if $Y$ is a
Banach space and $Z\subseteq Y^*$ is a norming subspace, then for
every strongly exposed point $y\in S_Y$, the exposing functional
belong to $Z$. Observe that this is the case of the ultraproduct
of the duals of a family of Banach space seen as a norm-closed
subspace of the dual of the ultraproduct of the spaces.

\begin{proof}
We only have to show that $(X_*)_\mathcal{U}$ has the Daugavet
property whenever $X_*$ does. Since $(X_*)_\mathcal{U}$ is the
predual of some $JBW^*$-triple, it suffices to show that its unit
ball has no strongly exposed points (Theorem~\ref{th-jbw}).
Therefore, we suppose, for the sake of contradiction, that the
unit ball of $(X_*)_\mathcal{U}$ has a strongly exposed point, say
$(x_i)$. By the preceding remark, there exists $(\phi_i)$ in the
unit sphere of $(X)_\mathcal{U}$ (which we can suppose to satisfy
$\|\phi_i\|=1$ for every $i$) which strongly expose $(x_i)$. Let
us fix $0<\varepsilon_0<1$. Now, for every $\alpha>0$, since $X_*$
has the Daugavet property, we can apply
Lemma~\ref{lemma-char}.$ii$ to get, for every $i\in I$, a point
$y_i\in S_{X_*}$ such that
$$\|x_i - y_i\|\geqslant 2 - \varepsilon_0 \qquad \text{and} \qquad
\re \phi_i(y_i)> 1 - \alpha/2.$$ Now, $(y_i)$ belong to the unit
ball of $(X_*)_\mathcal{U}$,
$$\|(x_i)-(y_i)\|= \lim_\mathcal{U}\|x_i-y_i\|\geqslant 2-\varepsilon_0,$$ and
\begin{equation*}
\re (\phi_i)[(y_i)]=\lim_\mathcal{U} \phi_i(y_i)> 1-\alpha.
\end{equation*}
Since $\alpha$ is arbitrary, we conclude that every slice of the
unit ball of $(X_*)_\mathcal{U}$ defined by $(\phi_i)$ has
diameter greater or equal than $2-\varepsilon_0$ (recall that $\re
(\phi_i)[(x_i)]=1$). Hence, $(\phi_i)$ does not strongly expose
$(x_i)$, a contradiction.
\end{proof}

As a consequence of the above theorem and Theorem~\ref{th-jbw}, we
obtain the following.

\begin{corollary}\label{cor-extremeultrafilter}
Let $X$ be a real or complex $JBW^*$-triple and $\mathcal{U}$ a
free ultrafilter on a set $I$. Then, the unit ball of
$(X_*)_\mathcal{U}$ have extreme points if and only if $B_{X_*}$
does.
\end{corollary}

As a consequence of Theorems \ref{th-jbw}, \ref{th-ultraproductJB}
and \ref{th-DaugtoultrapowerpredulaJBW}, we obtain

\begin{corollary}
Let $X_*$ be the predual of a real or complex $JBW^*$-triple $X$.
Then, $X_*$ has the uniform Daugavet property if and only if $X$
does.
\end{corollary}

It is worth mentioning that it is not known whether the uniform
Daugavet property passes from the dual of a Banach space to the
space.

\begin{remark}
The proof of Theorem~\ref{th-DaugtoultrapowerpredulaJBW} can be
straightforwardly adapted to show that {\slshape the ultraproduct
of a family of preduals of $JBW^*$-triples with the Daugavet
property also has the Daugavet property.\ } Therefore,
Corollary~\ref{cor-extremeultrafilter} can be also adapted to show
that {\slshape the unit ball of the ultraproduct of a family of
preduals of $JBW^*$-triples has no extreme points, provided that
the unit ball of each factor does not have any extreme point.}
\end{remark}

It is worth mentioning that Corollary~\ref{cor-extremeultrafilter}
can not be stated for general Banach spaces, as the following
example shows.

\begin{example}
{\slshape There exists a Banach space $X$ whose unit ball does not
have any extreme point and a free ultrafilter $\mathcal{U}$ on
$\mathbb{N}$ such that the unit ball of $X_\mathcal{U}$ has an
extreme point \cite[Example~2.14]{Hein2}.\ }
\end{example}

Let us comment a particular case in which the conclusion of
Corollary~\ref{cor-extremeultrafilter} can be easily stated.

\begin{remark}
{\slshape Let $X$ be a Banach space. Suppose that there exists
$\delta>0$ such that for every $x\in S_X$, there is $y\in X$ with
$\|y\|\geqslant \delta$ such that $\|x \pm y\|\leqslant 1$ (in
particular, $B_X$ has no extreme points). Then, for every free
ultrafilter $\mathcal{U}$ on a set $I$, the unit ball of
$X_\mathcal{U}$ does not have any extreme point.\ } Indeed, let
$(x_i)$ be a norm-one element of $X_\mathcal{U}$, which we can
suppose to satisfies $\|x_i\|=1$ for every $i$. Then, for every
$i\in I$, take $y_i\in X$ with $\|y_i\|\geqslant \delta$ and
$\|x_i \pm y_i\|\leqslant 1$. If we consider $(y_i)\in
X_\mathcal{U}$, then
$$\|(y_i)\|\geqslant \delta \qquad \text{and} \qquad \|(x_i) \pm
(y_i)\|\leqslant 1.$$ Therefore, $(x_i)$ is not an extreme point
of the unit ball of $X_\mathcal{U}$.
\end{remark}

It is easy to show that the above situation is fulfilled by
$L_1[0,1]$ with $\delta=1$.

\begin{example}
{\slshape For every $f\in L_1[0,1]$ with $\|f\|_1=1$, there is
$g\in L_1[0,1]$ with $\|g\|_1=1$ and such that $\|f \pm g\|_1=
1$.\ } Indeed, up to an isometric isomorphism, we can suppose
$f(t)\geqslant 0$ for every $t\in [0,1]$ and, by continuity, we
can find $t_0\in ]0,1[$ such that $$\int_0^{t_0} f(t)\ dt =
\int_{t_0}^1 f(t)\ dt = \dfrac{1}{2}.$$ Then, if we consider
$g=f\,(\chi_{[0,t_0]} - \chi_{[t_0,1]})\in L_1[0,1]$, we clearly
have $\|g\|_1=1$ and
\begin{align*}
\|f\pm g\|_1&=\int_0^{t_0} \bigl(f(t) \pm f(t)\bigr)\,dt\ + \
\int_{t_0}^1 \bigl(f(t) \mp f(t)\bigr)\,dt \\ &= \left(\dfrac12
\pm \dfrac12\right) + \left(\dfrac12 \mp \dfrac12\right)=1
\end{align*}
\end{example}

Actually, a very similar result (with $\delta$ arbitrarily closed
to $1$) can be stated for every $L_1(\mu)$ if $\mu$ does not have
any atom.

For the sake of completeness, we finish the paper by summarizing
the results of the present section in terms of $C^*$-algebras and
preduals of $W^*$-algebras.

\begin{corollary}$ $
\begin{enumerate}
\item[(a)] The ultraproduct of every family of real or complex
$C^*$-algebras with the Daugavet property also has the Daugavet
property. In particular, the Daugavet and the uniform Daugavet
property are equivalent for real or complex $C^*$-algebras.
\item[(b)] The ultrapower of a real or complex $C^*$-algebra has
atomic projections if and only if the algebra does.
\item[(c)] The ultraproduct of every family of preduals of real or
complex $W^*$-algebras with the Daugavet property also has the
Daugavet property. In particular, the Daugavet and the uniform
Daugavet property are equivalent for preduals of real or complex
$W^*$-algebras.
\item[(d)] Let $X_*$ be the predual of a real or complex
$W^*$-algebra $X$. Then, $X_*$ has the uniform Daugavet property
if and only if $X$ does.
\item[(e)] Let $Y$ be the predual of a real or complex
$W^*$-algebra. Then, $B_Y$ has an extreme point if and only if the
unit ball of every ultrapower of $Y$ does.
\end{enumerate}
\end{corollary}

\vspace{1cm}

\textbf{Acknowledgment:} The authors would like to express their
gratitude to \'{A}ngel Rodr\'{\i}guez Palacios for his valuable
suggestions, which have been decisive to the elaboration of this
paper. We also thank Gin\'{e}s L\'{o}pez, Javier Mer\'{\i}, Antonio Peralta,
and Armando Villena for helpful discussions. They also thanks the
referee, whose suggestions have improved the final version of this
paper.

\end{document}